\documentclass{amsart}
\usepackage{amssymb}
\usepackage[nobysame]{amsrefs}
\usepackage[usesnames,svgnames]{xcolor}
\usepackage[sc]{mathpazo}
\usepackage{tikz,pgf}
	\pgfdeclarelayer{background}
	\pgfsetlayers{background,main}
\usepackage{enumerate}
\usepackage{subcaption}
\usepackage{todonotes}
\usepackage[colorlinks=true,
	urlcolor=MidnightBlue,
	citecolor=DarkGreen]{hyperref}
\newtheorem{theorem}{Theorem}[section]
\newtheorem{proposition}[theorem]{Proposition}
\newtheorem{lemma}[theorem]{Lemma}
\newtheorem{corollary}[theorem]{Corollary}
\theoremstyle{remark}
\newtheorem{definition}[theorem]{Definition}
\newtheorem{example}[theorem]{Example}
\newtheorem{remark}[theorem]{Remark}

\newcommand{\defn}[1]{{\color{DarkGreen}\emph{#1}}}
\newcommand{\ie}{\text{i.e.}\;}
\newcommand{\CC}{\mathcal{C}}
\newcommand{\II}{\mathcal{I}}
\newcommand{\PP}{\mathcal{P}}
\newcommand{\QQ}{\mathcal{Q}}
\newcommand{\xx}{\mathbf{x}}
\newcommand{\yy}{\mathbf{y}}
\newcommand{\zz}{\mathbf{z}}
\DeclareMathOperator{\rk}{\mathrm rk}
\author{Henri M{\"u}hle}
\address{Institut f{\"u}r Algebra, Technische Universit{\"a}t Dresden, Zellescher Weg 12--14, 01069 Dresden, Germany.}
\email{henri.muehle@tu-dresden.de}
\thanks{This work was funded by a Public Grant overseen by the French National Research Agency (ANR) as part of the ``Investissements d'Avenir'' Program (Reference: ANR-10-LABX-0098).}
\title{On the Poset of Multichains}
\keywords{Chains, Multichains, Distributive Lattices, EL-Shellability, $P$-Partitions}
\subjclass[2010]{06A07 (primary), and 06B99 (secondary)}

\begin{document}

\begin{abstract}
	In this note we introduce the poset of $m$-multichains of a given poset $\mathcal{P}$.  Its elements are the multichains of $\mathcal{P}$ consisting of $m$ elements, and its partial order is the componentwise partial order of $\mathcal{P}$.  We show that this construction preserves a number of poset-theoretic and poset-topological properties of $\mathcal{P}$.  Moreover, we describe the structure of the poset of $m$-multichains of a finite distributive lattice, and provide a link to R.~Stanley's theory of $\mathcal{P}$-partitions.
\end{abstract}

\maketitle
\allowdisplaybreaks

\section{Introduction}
	\label{sec:introduction}
The incidence algebra of a locally finite poset $\PP=(P,\leq)$, as introduced by G.-C.~Rota~\cite{rota64foundations}, consists of all real-valued functions of two variables, say $f(x,y)$, where $x$ and $y$ range over the elements of $\PP$, and where $f(x,y)=0$ if $x\not\leq y$.  Among the elements of this algebra are two mutually inverse functions, the \defn{zeta function} and the \defn{M{\"o}bius function}, which are given by
\begin{displaymath}
	\zeta_{\PP}(x,y) = \begin{cases}1,&\text{if}\;x\leq y\\0,&\text{otherwise},\end{cases}\quad\text{and}\quad
	\mu_{\PP}(x,y)=\begin{cases}1,&\text{if}\;x=y\\-\sum\limits_{x\leq z<y}{\mu_{\PP}(x,z)},&\text{if}\;x<y\\0,&\text{otherwise},\end{cases}
\end{displaymath}
respectively.  Both of these functions are involved in the famous M{\"o}bius Inversion Formula~\cite{rota64foundations}*{Proposition~2}, which generalizes the Inclusion-Exclusion Principle.  Moreover, the M{\"o}bius function provides a deep link between combinatorics and algebraic topology: if $\PP$ has a least and a greatest element, then the value of the M{\"o}bius function given by these elements coincides with the reduced Euler characteristic of the order complex of the proper part of $\PP$~\cite{stanley11enumerative}*{Proposition~3.8.6}.  Much work has been done on the computation of the M{\"o}bius function for special posets, and many tools have been developed to aid this process.  We refer the reader to \cite{stanley11enumerative}*{Section~3} and \cites{bjorner96topological,wachs07poset}.

To some extent, the zeta function can be seen as a starting point for the work presented in this paper.  Suppose that $\PP$ is a poset with a least element $\hat{0}$ and greatest element $\hat{1}$.  Then $\zeta^{m}(\hat{0},\hat{1})$ is precisely the number of multichains of $\PP$ of length $m-1$.  This gives rise to the definition of the \defn{zeta polynomial} of $\PP$ by setting $Z(\PP,m)=\zeta^{m}(\hat{0},\hat{1})$.  Observe that this polynomial can be evaluated on all integers, since the zeta function is invertible.  See \cite{stanley11enumerative}*{Section~3.12} for more background and details.  Some more enumerative results involving the zeta polynomial of certain posets can be found in \cite{breckenridge93chains}.

We will, however, touch the enumerative aspect of the zeta polynomial only briefly, and rather focus on structural aspects of the set of $m$-multichains instead.  In particular, we view these multichains inside the $m^{\text{th}}$ direct power of $\PP$, which equips them naturally with a partial order; we denote this \defn{poset of $m$-multichains} by $\PP^{[m]}$.

A prominent example of a poset defined on the set of $m$-multichains of some other poset is D.~Armstrong's generalized noncrossing partition poset associated with a well-generated complex reflection group $W$~\cites{armstrong09generalized,bessis11cyclic}; usually denoted by $\mathcal{N\!C}_{W}^{(m)}$.  Its elements are the $m$-multichains of the noncrossing partition lattice $\mathcal{N\!C}_{W}$, but they are not ordered componentwise.  We remark that the poset $\mathcal{N\!C}_{W}^{[m]}$ of $m$-multichains of $\mathcal{N\!C}_{W}$ lives naturally inside the absolute order on the elements of the dual braid monoid of $W$~\cites{bessis03dual,bessis15finite}.  Other than for motivational reasons, we will not consider noncrossing partition posets in this note.  Another subposet of the $m^{\text{th}}$ direct power of $\PP$, which is also a subposet of $\PP^{[m]}$ is the \defn{$m$-cover poset} introduced and studied in \cite{kallipoliti15cover}.

We start in Section~\ref{sec:multichains} with the definition of the poset $\PP^{[m]}$, and we list first properties.  Subsequently, in Section~\ref{sec:quasi_varieties} we show that (quasi-)varieties of (finite) lattices are closed under the formation of the poset of $m$-multichains.  In Section~\ref{sec:topology} we investigate certain topological properties of $\PP^{[m]}$.  In particular, we prove that the existence of a certain well-behaved edge-labeling is preserved under the transition from $\PP$ to $\PP^{[m]}$.  We also show that the order complex of $\PP^{[m]}$ is contractible for $m>1$.  We conclude this note in Section~\ref{sec:applications} with an explicit description of the poset of $m$-multichains of a distributive lattice, and we relate this construction to the $\PP$-partitions of R.~Stanley.

\section{Structure of Multichain Posets}
	\label{sec:structure}
Let us recall some basic poset- and lattice-theoretic terminology.  For more background, we refer the reader to \cites{davey02introduction,gratzer11lattice}.  Throughout this note, we use the abbreviation $[n]=\{1,2,\ldots,n\}$ for some positive integer $n$.  All partially ordered sets (\defn{posets} for short) that we consider are supposed to be finite.  

\subsection{Posets of Multichains}
	\label{sec:multichains}
Let us start right away with the central definition of this note.

\begin{definition}\label{def:multichain_poset}
	Let $\PP=(P,\leq)$ be a poset.  For $m>0$ define the set of \defn{$m$-multichains} of $\PP$ by
	\begin{displaymath}
		P^{[m]} = \bigl\{(x_{1},x_{2},\ldots,x_{m})\mid x_{1}\leq x_{2}\leq\cdots\leq x_{m}\bigr\}.
	\end{displaymath}
	The \defn{poset of $m$-multichains} of $\PP$ is $\PP^{[m]}=\bigl(P^{[m]},\leq\bigr)$, where $\leq$ is considered componentwise, \ie $(x_{1},x_{2},\ldots,x_{m})\leq(y_{1},y_{2},\ldots,y_{m})$ if and only if $x_{i}\leq y_{i}$ for $i\in[m]$.
\end{definition}

Let $\PP=(P,\leq_{P})$ and $\QQ=(Q,\leq_{Q})$ be two posets.  Their \defn{direct product} is $\PP\times\QQ=(P\times Q,\leq)$, where $(x_{1},y_{1})\leq(x_{2},y_{2})$ if and only if $x_{1}\leq_{P}x_{2}$ and $y_{1}\leq_{Q}y_{2}$.  The direct product of $k$ posets $\PP_{1},\PP_{2},\ldots,\PP_{k}$ is abbreviated by $\prod\limits_{i=1}^{k}{\PP_{i}}$; if $\PP_{1}=\PP_{2}=\cdots=\PP_{k}=\PP$, then we write $\PP^{k}$ instead.

\begin{lemma}\label{lem:multichain_product}
	Let $\PP$ be a poset and let $m>0$.  The poset of $m$-multichains $\PP^{[m]}$ is an induced subposet of $\PP^{m}$.
\end{lemma}
\begin{proof}
	This is immediate from the definition.
\end{proof}

\begin{proposition}\label{prop:multichain_product}
	Let $\PP$ and $\QQ$ be two posets.  For $m>0$ we have $(\PP\times\QQ)^{[m]}\cong\PP^{[m]}\times\QQ^{[m]}$.
\end{proposition}
\begin{proof}
	Let $\PP=(P,\leq_{P})$ and $\QQ=(Q,\leq_{Q})$, and let $\bigl((x_{1},y_{1}),(x_{2},y_{2}),\ldots,(x_{m},y_{m})\bigr)\in(P\times Q)^{[m]}$.  By definition we have $x_{1}\leq_{P}x_{2}\leq_{P}\cdots\leq_{P}x_{m}$ and $y_{1}\leq_{Q}y_{2}\leq_{Q}\cdots\leq_{Q}y_{m}$, which implies that $\bigl((x_{1},x_{2},\ldots,x_{m}),(y_{1},y_{2},\ldots,y_{m})\bigr)\in P^{[m]}\times Q^{[m]}$.  It is clear that this is a bijection which we will denote by $f$.  It is also straightforward to verify that for any $\xx,\yy\in(P\times Q)^{[m]}$ we have $\xx\leq\yy$ if and only if $f(\xx)\leq f(\yy)$, which completes the proof.
\end{proof}

Two elements $x,y\in P$ form a \defn{cover relation} if $x<y$ and there is no $z\in P$ with $x<z<y$.  We usually write $x\lessdot y$ in this case.  We also say that $x$ is a \defn{lower cover} of $y$, or equivalently that $y$ is an \defn{upper cover} of $x$.

\begin{lemma}\label{lem:product_covers}
	Let $\PP=(P,\leq)$ be a poset, and let $m>0$.  Let $\xx,\yy\in P^{m}$, where $\xx=(x_{1},x_{2},\ldots,x_{m}),\yy=(y_{1},y_{2},\ldots,y_{m})$.  Then $\xx\lessdot\yy$ in $\PP^{m}$ if and only if there exists $j\in[m]$ with $x_{j}\lessdot y_{j}$ and $x_{i}=y_{i}$ for all $i\neq j$.
\end{lemma}
\begin{proof}
	Let $\xx,\yy\in P^{m}$ with $\xx=(x_{1},x_{2},\ldots,x_{m})$ and $\yy=(y_{1},y_{2},\ldots,y_{m})$, and suppose that $\xx<\yy$.  By construction there is at least one $j\in[m]$ with $x_{j}<y_{j}$.

	First assume that there is a unique $j\in[m]$ with $x_{j}<y_{j}$, and hence $x_{i}=y_{i}$ for $i\neq j$.  Let $\zz\in P^{[m]}$ with $\zz=(z_{1},z_{2},\ldots,z_{m})$.  If $\xx<\zz<\yy$, then by assumption we have $x_{i}=z_{i}=y_{i}$ for $i\neq j$, and therefore we have $x_{j}\leq z_{j}\leq y_{j}$.  If $x_{j}\lessdot y_{j}$, then we conclude $z_{j}=x_{j}$ or $z_{j}=y_{j}$, which implies $\xx\lessdot\yy$.  Otherwise there exists $z\in P$ with $x_{j}<z<y_{j}$ and by setting $z_{j}=z$, we see that $\xx$ and $\yy$ do not form a cover relation in $\PP^{[m]}$.
	
	Now suppose that there are at least two indices at which $\xx$ and $\yy$ differ.  Pick two of them, say $j_{1},j_{2}\in[m]$ with $j_{1}<j_{2}$.  We then have $x_{j_{1}}<y_{j_{1}}$ and $x_{j_{2}}<y_{j_{2}}$.  We can thus find elements $z\in P$ with $x_{j_{2}}\lessdot z\leq y_{j_{2}}$.  The $m$-tuple
	\begin{align*}
		\zz & = (x_{1},x_{2},\ldots,x_{j_{2}-1},z,x_{j_{2}+1},x_{j_{2}+2}\ldots,x_{m})
	\end{align*}
	certainly satisfies $\xx<\zz<\yy$.  We conclude that $\xx$ and $\yy$ do not form a cover relation in $\PP^{m}$.
\end{proof}

\begin{corollary}\label{cor:multichain_covers}
	Let $\PP=(P,\leq)$ be a poset, and let $m>0$.  Two elements $\xx,\yy\in\PP^{[m]}$ form a cover relation in $\PP^{[m]}$ if and only if they form a cover relation in $\PP^{m}$.
\end{corollary}
\begin{proof}
	The proof works more or less analogous to the proof of Lemma~\ref{lem:product_covers}.  We only need to be careful in the second case, when there are two indices $j_{1},j_{2}\in[m]$ with $x_{j_{1}}<y_{j_{1}}$ and $x_{j_{2}}<y_{j_{2}}$.  More precisely, we need to make sure that the element $\zz$ is indeed an $m$-multichain.  In fact, if $x_{j_{2}}<x_{j_{2}+1}$, then $\zz$ as defined above is in $P^{[m]}$.  Otherwise, if $x_{j_{2}}=x_{j_{2}+1}$, then by assumption 
	\begin{displaymath}
		x_{j_{2}+1}=x_{j_{2}}<y_{j_{2}}\leq y_{j_{2}+1}
	\end{displaymath}
	so that we can just start over with $j_{2}+1$ instead of $j_{2}$.  Once we hit the last index, we can construct $\zz$ without restrictions.
\end{proof}

An element $x\in P$ is \defn{minimal} if for every $y\in P$ with $y\leq x$ we have $x=y$.  Dually, $x$ is \defn{maximal} if for every $y\in P$ with $x\leq y$ we have $x=y$.  If $\PP$ has a unique minimal and a unique maximal element, then $\PP$ is \defn{bounded}.  In this case, we denote these elements by $\hat{0}$ and $\hat{1}$, respectively.

\begin{proposition}\label{prop:multichain_bounded}
	If $\PP$ is bounded, then $\PP^{[m]}$ is bounded for every $m>0$.
\end{proposition}
\begin{proof}
	This follows from the observation that for every minimal (maximal) element $x\in P$ the $m$-multichain $(x,x,\ldots,x)$ is minimal (maximal) in $\PP^{[m]}$.  
\end{proof}

Let $\xx\in P^{[m]}$.  If all entries of $\xx$ are distinct, then $\xx$ is an \defn{$m$-chain}, and its \defn{length} is $m-1$.  A chain is \defn{maximal} if it is maximal under inclusion.  A poset is \defn{graded} if all its maximal chains have the same length.

\begin{proposition}\label{prop:multichain_graded}
	If $\PP$ is graded, then $\PP^{[m]}$ is graded for every $m>0$.
\end{proposition}
\begin{proof}
	According to \cite{stanley11enumerative}*{Section~3.1} the assumption that $\PP$ is graded is equivalent to the existence of a rank function $\rk:P\to\mathbb{N}$ with $\rk(x)=0$ for every minimal element $x$, and $\rk(y)=\rk(x)+1$ whenever $x\lessdot y$.
	
	For $\xx\in P^{[m]}$ with $\xx=(x_{1},x_{2},\ldots,x_{m})$ define $\widehat{\rk}(\xx)=\sum\limits_{i=1}^{m}{\rk(x_{i})}$.  We claim that $\widehat{\rk}$ is a rank function of $\PP^{[m]}$.  
	
	Let $x\in P$ be minimal, and consider the minimal element $\xx=(x,x,\ldots,x)$ of $\PP^{[m]}$.  It follows immediately that $\widehat{\rk}(\xx)=0$.  Now let $\xx,\yy\in P^{[m]}$ with $\xx=(x_{1},x_{2},\ldots,x_{m})$ and $\yy=(y_{1},y_{2},\ldots,y_{m})$, and $\xx\lessdot\yy$.  Corollary~\ref{cor:multichain_covers} implies that there is a unique $j\in[m]$ with $x_{j}\lessdot y_{j}$ and $x_{i}=y_{i}$ for $i\neq j$.  We conclude
	\begin{displaymath}
		\widehat{\rk}(\yy) = \rk(y_{j}) + \sum_{i\neq j}{\rk(y_{i})} = \rk(x_{j})+1 + \sum_{i\neq j}{\rk(x_{i})} = \widehat{\rk}(\xx)+1
	\end{displaymath}
	as desired.
\end{proof}

A \defn{lattice} is a poset $\PP=(P,\leq)$ in which every two elements $x,y\in P$ have a greatest lower bound (their \defn{meet}; denoted by $x\wedge y$) and a least upper bound (their \defn{join}; denoted by $x\vee y$).  Let $\PP=(P,\leq_{P})$ and $\QQ=(Q,\leq_{Q})$ be two lattices.  We say that $\PP$ is a \defn{sublattice} of $\QQ$ if $P\subseteq Q$ and for every $x,y\in P$ we have $x\wedge_{P}y=x\wedge_{Q}y$ and $x\vee_{P}y=x\vee_{Q}y$.  

\begin{theorem}\label{thm:multichain_sublattice_product}
	If $\PP$ is a lattice, then $\PP^{[m]}$ is a sublattice of $\PP^{m}$ for every $m>0$. 
\end{theorem}
\begin{proof}
	Let $\PP=(P,\leq)$ be a lattice.  Lemma~\ref{lem:multichain_product} states that $\PP^{[m]}$ is a subposet of $\PP^{m}$.  It remains to show that meets and joins are preserved.  
	
	Let $\xx,\yy\in P^{[m]}$ with $\xx=(x_{1},x_{2},\ldots,x_{m})$ and $\yy=(y_{1},y_{2},\ldots,y_{m})$.  For $i\in[m-1]$ we have $x_{i}\leq x_{i+1}$ and $y_{i}\leq y_{i+1}$, which implies $x_{i}\wedge y_{i}\leq x_{i+1}\wedge y_{i+1}$ and $x_{i}\vee y_{i}\leq x_{i+1}\vee y_{i+1}$.  We conclude that the componentwise meet (join) of $\xx$ and $\yy$ (which is the meet (join) of $\xx$ and $\yy$ in $\PP^{m}$) is contained in $\PP^{[m]}$, and it must thus be the meet (join) of $\xx$ and $\yy$ in $\PP^{[m]}$.
\end{proof}

\begin{corollary}\label{cor:multichain_sublattice}
	If $\PP$ is a sublattice of $\QQ$, then $\PP^{[m]}$ is a sublattice of $\QQ^{[m]}$ for every $m>0$.
\end{corollary}
\begin{proof}
	This follows directly from Theorem~\ref{thm:multichain_sublattice_product}.
\end{proof}

Let $\PP=(P,\leq_{P})$ and $\QQ=(Q,\leq_{Q})$ be two lattices.  A map $f:P\to Q$ is a \defn{lattice homomorphism} if $f(x\wedge_{P}y)=f(x)\wedge_{Q}f(y)$ and $f(x\vee_{P}y)=f(x)\vee_{Q}f(y)$ for every $x,y\in P$.  If $f$ is a surjective lattice homomorphism, then $\QQ$ is a \defn{homomorphic image} of $\PP$.

\begin{corollary}\label{cor:multichain_homomorphic_image}
	If $\QQ$ is a homomorphic image of $\PP$, then $\QQ^{[m]}$ is a homomorphic image of $\PP^{[m]}$ for every $m>0$. 
\end{corollary}
\begin{proof}
	Let $\PP=(P,\leq_{P})$ and $\QQ=(Q,\leq_{Q})$ be two lattices, and let $f:P\to Q$ be a surjective lattice homomorphism.  Define $\hat{f}:P^{[m]}\to Q^{[m]}$ componentwise, \ie
	\begin{displaymath}
		\hat{f}\bigl((x_{1},x_{2},\ldots,x_{m})\bigr) = \bigl(f(x_{1}),f(x_{2}),\ldots,f(x_{m})\bigr).
	\end{displaymath}
	This is well defined, since for $i,j\in[m]$ with $i\leq j$ we have $x_{i}\leq_{P}x_{j}$, which is equivalent to $x_{j}=x_{i}\vee_{P}x_{j}$.  We thus have $f(x_{j})=f(x_{i}\vee_{P}x_{j})=f(x_{i})\vee_{Q}f(x_{j})$, which is equivalent to $f(x_{i})\leq_{Q}f(x_{j})$.  Since $f$ is surjective, so is $\hat{f}$. 
	
	Let $\xx,\yy\in P^{[m]}$ with $\xx=(x_{1},x_{2},\ldots,x_{m})$ and $\yy=(y_{1},y_{2},\ldots,y_{m})$.  Theorem~\ref{thm:multichain_sublattice_product} implies together with the assumption that $f$ is a lattice homomorphism that
	\begin{align*}
		\hat{f}(\xx\wedge_{P}\yy) & = \bigl(f(x_{1}\wedge_{P}y_{1}),f(x_{2}\wedge_{P}y_{2}),\ldots,f(x_{m}\wedge_{P}y_{m})\bigr)\\
		& = \bigl(f(x_{1})\wedge_{Q}f(y_{1}),f(x_{2})\wedge_{Q}f(y_{2}),\ldots,f(x_{m})\wedge_{Q}f(y_{m})\bigr)\\
		& = \hat{f}(\xx)\wedge_{Q}\hat{f}(\yy),
	\end{align*}
	and likewise for joins.
\end{proof}

Let $\PP=(P,\leq)$ be a lattice.  An equivalence relation $\Theta$ on $P$ is a \defn{lattice congruence} if for every $x,y,z\in P$ we have $(x,y)\in\Theta$ implies $(x\wedge z,y\wedge z)\in\Theta$ and $(x\vee z,y\vee z)\in\Theta$.  The equivalence classes of $\Theta$ together with the induced order form the \defn{quotient lattice} $\PP/\Theta$.  

If $\Theta$ is a lattice congruence of $\PP=(P,\leq)$, then the map $f:P\to P/\Theta, x\mapsto [x]_{\Theta}$ is a surjective lattice homomorphism.  Conversely, for lattices $\PP=(P,\leq_{P})$ and $\QQ=(Q,\leq_{Q})$ any surjective lattice homomorphism $f:P\to Q$ induces a lattice congruence $\Theta$ such that $\PP/\Theta\cong\QQ$.  The equivalence classes of $\Theta$ are precisely the fibers (preimages) of $f$.

\begin{corollary}\label{cor:multichain_quotient}
	If $\QQ$ is a quotient lattice of $\PP$, then $\QQ^{[m]}$ is a quotient lattice of $\PP^{[m]}$ for every $m>0$.
\end{corollary}
\begin{proof}
	This follows from Corollary~\ref{cor:multichain_homomorphic_image}.
\end{proof}

Given a poset $\PP=(P,\leq_{P})$ its \defn{dual} is the poset $\PP^{d}=(P,\preceq_{P})$ with $x\preceq_{P}y$ if and only if $y\leq_{P}x$.  We say that a poset $\QQ=(Q,\leq_{Q})$ is the \defn{dual} of $\PP$ if $\QQ\cong\PP^{d}$.  In particular there exists a bijective lattice homomorphism $f:P\to Q$ such that $x\leq_{P}y$ if and only if $f(y)\leq_{Q}f(x)$.

\begin{corollary}\label{cor:multichain_dual}
	If $\PP$ and $\QQ$ are dual posets, then $\PP^{[m]}$ and $\QQ^{[m]}$ are dual posets for every $m>0$.
\end{corollary}
\begin{proof}
	This follows from Corollary~\ref{cor:multichain_homomorphic_image}.
\end{proof}

\subsection{(Quasi-)Varieties of Lattices}
	\label{sec:quasi_varieties}
Now fix a lattice $\PP=(P,\leq)$.  Let us now recursively define \defn{lattice terms}.  Any $x\in P$ is a lattice term of length $1$, and if $t_{1},t_{2},\ldots,t_{s}$ are lattice terms of lengths $k_{1},k_{2},\ldots,k_{s}$, then $(t_{1}\wedge t_{2}\wedge\cdots\wedge t_{s})$ and $(t_{1}\vee t_{2}\vee\cdots\vee t_{s})$ are both lattice terms of lengths $1+k_{1}+k_{2}+\cdots+k_{s}$.  A \defn{quasi-identity} is an implication of the form
\begin{displaymath}
	s_{1}=t_{1}\;\text{and}\;s_{2}=t_{2}\;\text{and}\;\ldots\;\text{and}\;s_{n}=t_{n}\quad\text{imply}\quad s=t,
\end{displaymath}
where $s_{1},s_{2},\ldots,s_{n},t_{1},t_{2},\ldots,t_{n}$ and $s,t$ are all lattice terms.  If $n=0$, then we simply speak of an \defn{identity}.  A class $\mathbf{K}$ of lattices is a \defn{(quasi-)variety} if its members can be completely described by a set of (quasi-)identities.

\begin{example}\label{ex:lattices}
	Lattices can be characterized equivalently as a set with two binary operations $\wedge$ and $\vee$ that each are associative, commutative, and idempotent, and that additionally satisfy the absorption laws, see for instance~\cite{gratzer11lattice}*{Section~1.10}.  These ``lattice axioms'' are identities, which implies that the class of all lattices is a variety.
\end{example}

\begin{example}\label{ex:distributive}
	A lattice $\PP=(P,\leq)$ is \defn{distributive} if $x\wedge(y\vee z)=(x\wedge y)\vee(x\wedge z)$ and $x\vee(y\wedge z)=(x\vee y)\wedge(x\vee z)$ holds for all $x,y,z\in P$.  Therefore, the class of all distributive lattices is a variety.
\end{example}

\begin{example}\label{ex:modular}
	A lattice $\PP=(P,\leq)$ is \defn{modular} if $(x\wedge z)\vee(y\wedge z)=\bigl((x\wedge z)\vee y\bigr)\wedge z$ holds for all $x,y,z\in P$.  Therefore, the class of all modular lattices is a variety.
\end{example}

\begin{example}\label{ex:semidistributive}
	A lattice $\PP=(P,\leq)$ is \defn{join-semidistributive} if for all $x,y,z\in P$ the following implication is satisfied: if $x\vee y=x\vee z$, then $x\vee(y\wedge z)=x\vee y$.  Therefore, the class of all join-semidistributive lattices is a quasi-variety.  We can dually define the quasi-variety of meet-semidistributive lattices.
\end{example}

\begin{example}\label{ex:semimodular}
	A lattice $\PP=(P,\leq)$ is \defn{semimodular} if for all $x,y\in P$ the following implication is satisfied: if $x\wedge y\lessdot x$, then $y\lessdot x\vee y$.  Observe that this is not a quasi-identity, since it involves cover relations.  
	
	Since we are dealing only with finite lattices, however, we can find another characterization of semimodularity in terms of quasi-identities.  We say that $\PP$ is \defn{$M$-symmetric} if for all $x,y,z\in P$ the following implications are satisfied: 
	\begin{align*}
		x\leq z\quad\text{implies}\quad x\vee(y\wedge z)=(x\vee y)\wedge z,\\
		x\leq y\quad\text{implies}\quad x\vee(y\wedge z)=(x\vee z)\wedge y.
	\end{align*}
	Since $x\leq z$ and $x\leq y$ are equivalent to $x\wedge z=x$ and $x\wedge y=x$, respectively, we conclude that these implications are quasi-identities.  Therefore, the class of $M$-symmetric lattices is a quasi-variety.  For finite lattices, the notions of semimodularity and $M$-symmetry agree~\cite{malliah86equivalence}*{Theorem~3.1}.  Therefore, the class of finite semimodular lattices is a quasi-variety.
\end{example}

\begin{theorem}[\cite{birkhoff35on}*{Theorems~6~and~7}]\label{thm:variety_characterization}
	A class $\mathbf{K}$ of lattices is a variety if and only if it is closed under the formation of homomorphic images, sublattices, and direct products.
\end{theorem}

\begin{theorem}[\cite{burris81course}*{Theorem~2.25}]\label{thm:quasivariety_characterization}
	A class $\mathbf{K}$ of finite lattices is a quasi-variety if and only if it is closed under the formation of sublattices, and direct products.
\end{theorem}

We obtain the following result.

\begin{theorem}\label{thm:multichain_quasivariety}
	Let $\mathbf{K}$ be a class of finite lattices.  If $\mathbf{K}$ is a (quasi-)variety and $\PP\in\mathbf{K}$, then $\PP^{[m]}\in\mathbf{K}$ for every $m>0$.
\end{theorem}
\begin{proof}
	This follows from Theorems~\ref{thm:variety_characterization} and \ref{thm:quasivariety_characterization} in conjunction with Theorem~\ref{thm:multichain_sublattice_product}.
\end{proof}

\begin{remark}
	The notion of a (quasi-)variety exists in fact on the level of universal algebra, where it can be defined for an arbitrary algebraic structure.  Since this note deals exclusively with posets and lattices, however, we decided to present the content of this section tailored to this particular situation.
\end{remark}

\section{Topology of Multichain Posets}
	\label{sec:topology}
There is a natural way to associate a topological space with a poset $\PP$ via the geometric realization of the \defn{order complex} of $\PP$, \ie the simplicial complex whose faces are the chains of $\PP$.  Poset topology is the mathematical discipline that studies topological properties of this simplicial complex from the poset perspective.  See for instance \cites{bjorner96topological,wachs07poset} for an introduction to this topic.

Important tools in the study of the topology of the order complex of a poset are certain edge-labelings, which have their origin in \cites{stanley74finite,bjorner80shellable}.  See also \cites{bjorner96shellable,bjorner97shellable} for further background.  

Let $\PP=(P,\leq)$ be a bounded poset, and let $\mathcal{E}(\PP)=\{(x,y)\mid x\lessdot y\}$ denote its set of cover relations.  An \defn{edge-labeling} is simply a map $\lambda:\mathcal{E}(\PP)\to\Lambda$, where $(\Lambda,\leq_{\Lambda})$ is an arbitrary poset.  Given a maximal chain $\xx=(x_{1},x_{2},\ldots,x_{m})$ of $\PP$, its label sequence is $\lambda(\xx)=\bigl(\lambda(x_{1},x_{2}),\lambda(x_{1},x_{2}),\ldots,\lambda(x_{m-1},x_{m})\bigr)$.  We say that $\xx$ is \defn{rising} if $\lambda(\xx)$ is strictly increasing, and $\xx$ is \defn{falling} if $\lambda(\xx)$ is weakly decreasing.  

The set of tuples over $\Lambda$ is defined by
\begin{displaymath}
	\Lambda^{*}=\bigl\{(l_{1},l_{2},\ldots,l_{s})\mid s\in\mathbb{N},l_{i}\in\Lambda\;\text{for}\;i\in[s]\bigr\}=\bigcup_{s\geq 0}{\Lambda^{s}}.
\end{displaymath}
We consider $\Lambda^{*}$ equipped with the lexicographic order $\leq_{\text{lex}}$ which is defined by
\begin{displaymath}
	(k_{1},k_{2},\ldots,k_{s})\leq_{\text{lex}}(l_{1},l_{2},\ldots,l_{t})
\end{displaymath}
if and only if either $s\leq t$ and $k_{i}=l_{i}$ for $i\in[s]$, or $k_{i}<_{\Lambda}l_{i}$ for the least $i$ with $k_{i}\neq l_{i}$.  For two maximal chains $\xx,\xx'$ of $\PP$ we say that $\xx$ \defn{precedes} $\xx'$ if $\lambda(\xx)\leq_{\text{lex}}\lambda(\xx')$.

An \defn{interval} of $\PP$ is a set of the form $[x,y]=\{z\in P\mid x\leq z\leq y\}$ for $x,y\in P$.  An edge-labeling $\lambda$ of $\PP$ is an \defn{EL-labeling} if for every interval $[x,y]$ there exists a unique rising maximal chain in $[x,y]$, which precedes every other maximal chain of $[x,y]$.  A bounded poset that admits an EL-labeling is \defn{EL-shellable}.

Let us briefly outline the importance of EL-labelings.  The \defn{proper part} of a bounded poset $\PP=(P,\leq)$ is the induced subposet $\bar{\PP}=\bigl(P\setminus\{\hat{0},\hat{1}\},\leq\bigr)$.

\begin{theorem}[\cite{bjorner96shellable}*{Theorems~5.8~and~5.9}]\label{thm:el_shellable}
	If $\PP$ is a EL-shellable poset, then the order complex of $\bar{\PP}$ is shellable.  The $i^{\text{th}}$ Betti number of this order complex is given by the number of falling maximal chains of $\PP$ of length $i+2$.
\end{theorem}

\begin{theorem}[\cite{bjorner96shellable}*{Proposition~5.7}]\label{thm:mobius}
	If $\PP$ is a EL-shellable poset, then $\mu_{\PP}(\hat{0},\hat{1})$ equals the number of falling maximal chains of $\PP$ of even length minus the number of falling maximal chains of $\PP$ of odd length.
\end{theorem}

Suppose that $\PP_{1},\PP_{2},\ldots,\PP_{m}$ are EL-shellable posets such that $\lambda_{i}:\mathcal{E}(\PP_{i})\to\Lambda_{i}$ is an EL-labeling of $\PP_{i}$ for $i\in[m]$.  Fix a formal symbol $\delta$, and adjoin it as a least element to $\Lambda_{i}$ for $i\in[m]$.  The \defn{product labeling} $\lambda^{m}:\mathcal{E}\Bigl(\prod\limits_{i=1}^{m}{\PP_{i}}\Bigr)\to\prod\limits_{i=1}^{m}{\bigl(\Lambda_{i}\cup\{\delta\}\bigr)}$ is defined by
\begin{displaymath}
	\lambda^{m}\bigl((x_{1},x_{2},\ldots,x_{m}),(y_{1},y_{2},\ldots,y_{m})\bigr) = \bigl(\delta,\delta,\ldots,\delta,\lambda_{j}(x_{j},y_{j}),\delta,\delta,\ldots,\delta\bigr),
\end{displaymath}
where $j$ is the unique index in $[m]$ with $x_{j}\lessdot y_{j}$ from Lemma~\ref{lem:product_covers}.

\begin{theorem}[\cite{bjorner97shellable}*{Proposition~10.15}]\label{thm:product_shellable}
	Let $\PP_{1},\PP_{2},\ldots,\PP_{m}$ be EL-shellable posets with EL-labelings $\lambda_{1},\lambda_{2},\ldots,\lambda_{m}$.  The product labeling $\lambda^{m}$ is an EL-labeling of $\prod\limits_{i=1}^{m}{\PP_{i}}$.
\end{theorem}

We have the following result.

\begin{theorem}\label{thm:multichain_shellable}
	Let $\PP$ be an EL-shellable poset with EL-labeling $\lambda$.  The product labeling $\lambda^{m}$ is an EL-labeling of $\PP^{[m]}$ for $m>0$.
\end{theorem}
\begin{proof}
	Let $\PP=(P,\leq)$ and fix $\xx,\yy\in P^{[m]}$ with $\xx=(x_{1},x_{2},\ldots,x_{m})$ and $\yy=(y_{1},y_{2},\ldots,y_{m})$ and $\xx\leq\yy$.  For $i\in[m]$ define
	\begin{displaymath}
		\zz_{i} = (x_{1},x_{2},\ldots,x_{m-i},y_{m-i+1},y_{m-i+2},\ldots,y_{m})\in P^{[m]}.
	\end{displaymath}
	Let $C_{i}$ denote the unique rising chain in the interval $[x_{m-i+1},y_{m-i+1}]$ of $\PP$.  The concatenation $C$ of the chains $C_{1},C_{2},\ldots,C_{m}$ is thus a maximal chain in the interval $[\xx,\yy]$ of $\PP^{[m]}$ which is rising with respect to $\lambda^{m}$.  Moreover, $C$ passes through $\zz_{1},\zz_{2},\ldots,\zz_{m}$, and since $\lambda$ is an EL-labeling of $\PP$ no other maximal chain of $[\xx,\yy]$ passing through $\zz_{1},\zz_{2},\ldots,\zz_{m}$ can be rising.  Moreover, $C$ precedes any maximal chain of $[\xx,\yy]$ passing through $\zz_{1},\zz_{2},\ldots,\zz_{m}$.
	
	Let $C'$ be another maximal chain in $\PP^{[m]}$, and let $i$ be the first index such that $C'$ does not pass through $\zz_{i}$.  There must thus be elements $\xx_{1},\xx_{2},\xx_{3},\xx_{4}$ in $C'$ with $\xx_{1}\lessdot\xx_{2}\leq\xx_{3}\lessdot\xx_{4}$ such that $\lambda^{m}(\xx_{1},\xx_{2})$ has an entry $k\neq\delta$ at position $j<m-i$, and $\lambda^{m}(\xx_{3},\xx_{4})$ has an entry $l\neq\delta$ at position $m-i$.  Therefore $C'$ cannot be rising, and it cannot precede $C$.
	
	We conclude that $\lambda^{m}$ is an EL-labeling of $\PP^{[m]}$.
\end{proof}

See Figures~\ref{fig:poset_example_1} and \ref{fig:poset_example_2} for illustrations of Theorem~\ref{thm:multichain_shellable}. 

\begin{figure}
	\centering
	\begin{subfigure}[t]{.3\textwidth}
		\centering
		\begin{tikzpicture}\small
			\def\x{1.5};
			\def\y{1.5};
			\draw(2*\x,1*\y) node(n1){$\hat{0}$};
			\draw(1*\x,2*\y) node(n2){$a$};
			\draw(2*\x,2*\y) node(n3){$b$};
			\draw(3*\x,2*\y) node(n4){$c$};
			\draw(2*\x,3*\y) node(n5){$\hat{1}$};
			\draw(n1) -- (n2) node[fill=white,inner sep=.8] at(1.5*\x,1.5*\y){\tiny\color{white!50!black}$1$};
			\draw(n1) -- (n3) node[fill=white,inner sep=.8] at(2*\x,1.5*\y){\tiny\color{white!50!black}$2$};
			\draw(n1) -- (n4) node[fill=white,inner sep=.8] at(2.5*\x,1.5*\y){\tiny\color{white!50!black}$3$};
			\draw(n2) -- (n5) node[fill=white,inner sep=.8] at(1.5*\x,2.5*\y){\tiny\color{white!50!black}$3$};
			\draw(n3) -- (n5) node[fill=white,inner sep=.8] at(2*\x,2.5*\y){\tiny\color{white!50!black}$1$};
			\draw(n4) -- (n5) node[fill=white,inner sep=.8] at(2.5*\x,2.5*\y){\tiny\color{white!50!black}$2$};
		\end{tikzpicture}
		\caption{A bounded poset.}
		\label{fig:p}
	\end{subfigure}
	\hspace*{.5cm}
	\begin{subfigure}[t]{.6\textwidth}
		\centering
		\begin{tikzpicture}\small
			\def\x{1.5};
			\def\y{1.5};
			\draw(3*\x,1*\y) node(n1){$(\hat{0},\hat{0})$};
			\draw(2*\x,2*\y) node(n2){$(\hat{0},a)$};
			\draw(3*\x,2*\y) node(n3){$(\hat{0},b)$};
			\draw(4*\x,2*\y) node(n4){$(\hat{0},c)$};
			\draw(1*\x,3*\y) node(n5){$(a,a)$};
			\draw(2*\x,3*\y) node(n6){$(\hat{0},\hat{1})$};
			\draw(4*\x,3*\y) node(n7){$(b,b)$};
			\draw(5*\x,3*\y) node(n8){$(c,c)$};
			\draw(2*\x,4*\y) node(n9){$(a,\hat{1})$};
			\draw(3*\x,4*\y) node(n10){$(b,\hat{1})$};
			\draw(4*\x,4*\y) node(n11){$(c,\hat{1})$};
			\draw(3*\x,5*\y) node(n12){$(\hat{1},\hat{1})$};
			\draw(n1) -- (n2) node[fill=white,inner sep=.8] at(2.5*\x,1.5*\y){\tiny\color{white!50!black}$(\delta,1)$};
			\draw(n1) -- (n3) node[fill=white,inner sep=.8] at(3*\x,1.5*\y){\tiny\color{white!50!black}$(\delta,2)$};
			\draw(n1) -- (n4) node[fill=white,inner sep=.8] at(3.5*\x,1.5*\y){\tiny\color{white!50!black}$(\delta,3)$};
			\draw(n2) -- (n5) node[fill=white,inner sep=.8] at(1.5*\x,2.5*\y){\tiny\color{white!50!black}$(1,\delta)$};
			\draw(n2) -- (n6) node[fill=white,inner sep=.8] at(2*\x,2.5*\y){\tiny\color{white!50!black}$(\delta,3)$};
			\draw(n3) -- (n6) node[fill=white,inner sep=.8] at(2.5*\x,2.5*\y){\tiny\color{white!50!black}$(\delta,1)$};
			\draw(n3) -- (n7) node[fill=white,inner sep=.8] at(3.5*\x,2.5*\y){\tiny\color{white!50!black}$(2,\delta)$};
			\draw(n4) -- (n6) node[fill=white,inner sep=.8] at(3*\x,2.5*\y){\tiny\color{white!50!black}$(\delta,2)$};
			\draw(n4) -- (n8) node[fill=white,inner sep=.8] at(4.5*\x,2.5*\y){\tiny\color{white!50!black}$(3,\delta)$};
			\draw(n5) -- (n9) node[fill=white,inner sep=.8] at(1.5*\x,3.5*\y){\tiny\color{white!50!black}$(\delta,3)$};
			\draw(n6) -- (n9) node[fill=white,inner sep=.8] at(2*\x,3.5*\y){\tiny\color{white!50!black}$(1,\delta)$};
			\draw(n6) -- (n10) node[fill=white,inner sep=.8] at(2.5*\x,3.5*\y){\tiny\color{white!50!black}$(2,\delta)$};
			\draw(n6) -- (n11) node[fill=white,inner sep=.8] at(3*\x,3.5*\y){\tiny\color{white!50!black}$(3,\delta)$};
			\draw(n7) -- (n10) node[fill=white,inner sep=.8] at(3.5*\x,3.5*\y){\tiny\color{white!50!black}$(\delta,1)$};
			\draw(n8) -- (n11) node[fill=white,inner sep=.8] at(4.5*\x,3.5*\y){\tiny\color{white!50!black}$(\delta,2)$};
			\draw(n9) -- (n12) node[fill=white,inner sep=.8] at(2.5*\x,4.5*\y){\tiny\color{white!50!black}$(3,\delta)$};
			\draw(n10) -- (n12) node[fill=white,inner sep=.8] at(3*\x,4.5*\y){\tiny\color{white!50!black}$(1,\delta)$};
			\draw(n11) -- (n12) node[fill=white,inner sep=.8] at(3.5*\x,4.5*\y){\tiny\color{white!50!black}$(2,\delta)$};
		\end{tikzpicture}
		\caption{The poset of $2$-multichains of the poset in Figure~\ref{fig:p}.}
		\label{fig:p2}
	\end{subfigure}
	\caption{An EL-shellable poset and its poset of $2$-multichains.}
	\label{fig:poset_example_1}
\end{figure}

\begin{figure}
	\centering
	\begin{subfigure}[t]{.31\textwidth}
		\centering
		\begin{tikzpicture}\small
			\def\x{1.1};
			\def\y{1.1};
			\draw(2*\x,1*\y) node(n1){$\hat{0}$};
			\draw(1*\x,2*\y) node(n2){$a$};
			\draw(3*\x,2.5*\y) node(n3){$b$};
			\draw(1*\x,3*\y) node(n4){$c$};
			\draw(2*\x,4*\y) node(n5){$\hat{1}$};
			\draw(n1) -- (n2) node[fill=white,inner sep=.8] at(1.5*\x,1.5*\y){\tiny\color{white!50!black}$1$};
			\draw(n1) -- (n3) node[fill=white,inner sep=.8] at(2.5*\x,1.75*\y){\tiny\color{white!50!black}$3$};
			\draw(n2) -- (n4) node[fill=white,inner sep=.8] at(1*\x,2.5*\y){\tiny\color{white!50!black}$2$};
			\draw(n3) -- (n5) node[fill=white,inner sep=.8] at(2.5*\x,3.25*\y){\tiny\color{white!50!black}$1$};
			\draw(n4) -- (n5) node[fill=white,inner sep=.8] at(1.5*\x,3.5*\y){\tiny\color{white!50!black}$3$};
		\end{tikzpicture}
		\caption{Another bounded poset.}
		\label{fig:q}
	\end{subfigure}
	\hspace*{.5cm}
	\begin{subfigure}[t]{.61\textwidth}
		\centering
		\begin{tikzpicture}\small
			\def\x{1.1};
			\def\y{1.1};
			\draw(3*\x,1*\y) node(n1){$(\hat{0},\hat{0})$};
			\draw(2*\x,2*\y) node(n2){$(\hat{0},a)$};
			\draw(4*\x,2.5*\y) node(n3){$(\hat{0},b)$};
			\draw(1*\x,3*\y) node(n4){$(a,a)$};
			\draw(2*\x,3*\y) node(n5){$(\hat{0},c)$};
			\draw(1*\x,4*\y) node(n6){$(a,c)$};
			\draw(3*\x,4*\y) node(n7){$(\hat{0},\hat{1})$};
			\draw(4*\x,4*\y) node(n8){$(b,b)$};
			\draw(1*\x,5*\y) node(n9){$(c,c)$};
			\draw(2*\x,5*\y) node(n10){$(a,\hat{1})$};
			\draw(4*\x,5.5*\y) node(n11){$(b,\hat{1})$};
			\draw(2*\x,6*\y) node(n12){$(c,\hat{1})$};
			\draw(3*\x,7*\y) node(n13){$(\hat{1},\hat{1})$};
			\draw(n1) -- (n2) node[fill=white, inner sep=.8] at(2.5*\x,1.5*\y){\tiny\color{white!50!black}$(\delta,1)$};
			\draw(n1) -- (n3) node[fill=white, inner sep=.8] at(3.5*\x,1.75*\y){\tiny\color{white!50!black}$(\delta,3)$};
			\draw(n2) -- (n4) node[fill=white, inner sep=.8] at(1.5*\x,2.5*\y){\tiny\color{white!50!black}$(1,\delta)$};
			\draw(n2) -- (n5) node[fill=white, inner sep=.8] at(2*\x,2.5*\y){\tiny\color{white!50!black}$(\delta,2)$};
			\draw(n3) -- (n7) node[fill=white, inner sep=.8] at(3.5*\x,3.25*\y){\tiny\color{white!50!black}$(\delta,1)$};
			\draw(n3) -- (n8) node[fill=white, inner sep=.8] at(4*\x,3.25*\y){\tiny\color{white!50!black}$(3,\delta)$};
			\draw(n4) -- (n6) node[fill=white, inner sep=.8] at(1*\x,3.5*\y){\tiny\color{white!50!black}$(\delta,2)$};
			\draw(n5) -- (n6) node[fill=white, inner sep=.8] at(1.5*\x,3.5*\y){\tiny\color{white!50!black}$(1,\delta)$};
			\draw(n5) -- (n7) node[fill=white, inner sep=.8] at(2.5*\x,3.5*\y){\tiny\color{white!50!black}$(\delta,3)$};
			\draw(n6) -- (n9) node[fill=white, inner sep=.8] at(1*\x,4.5*\y){\tiny\color{white!50!black}$(2,\delta)$};
			\draw(n6) -- (n10) node[fill=white, inner sep=.8] at(1.5*\x,4.5*\y){\tiny\color{white!50!black}$(\delta,3)$};
			\draw(n7) -- (n10) node[fill=white, inner sep=.8] at(2.5*\x,4.5*\y){\tiny\color{white!50!black}$(1,\delta)$};
			\draw(n7) -- (n11) node[fill=white, inner sep=.8] at(3.5*\x,4.75*\y){\tiny\color{white!50!black}$(3,\delta)$};
			\draw(n8) -- (n11) node[fill=white, inner sep=.8] at(4*\x,4.75*\y){\tiny\color{white!50!black}$(\delta,1)$};
			\draw(n9) -- (n12) node[fill=white, inner sep=.8] at(1.5*\x,5.5*\y){\tiny\color{white!50!black}$(\delta,3)$};
			\draw(n10) -- (n12) node[fill=white, inner sep=.8] at(2*\x,5.5*\y){\tiny\color{white!50!black}$(2,\delta)$};
			\draw(n11) -- (n13) node[fill=white, inner sep=.8] at(3.5*\x,6.25*\y){\tiny\color{white!50!black}$(1,\delta)$};
			\draw(n12) -- (n13) node[fill=white, inner sep=.8] at(2.5*\x,6.5*\y){\tiny\color{white!50!black}$(3,\delta)$};
		\end{tikzpicture}
		\caption{The poset of $2$-multichains of the poset in Figure~\ref{fig:q}.}
		\label{fig:q2}
	\end{subfigure}
	\caption{Another EL-shellable poset and its poset of $2$-multichains.}
	\label{fig:poset_example_2}
\end{figure}

Let $\PP=(P,\leq)$ be bounded.  An element $x\in P$ is an \defn{atom} if $\hat{0}\lessdot x$ and it is a \defn{coatom} if $x\lessdot\hat{1}$.  For $m>0$, define $\mathbf{\hat{0}}=(\hat{0},\hat{0},\ldots,\hat{0})\in P^{[m]}$ and $\mathbf{\hat{1}}=(\hat{1},\hat{1},\ldots,\hat{1})\in P^{[m]}$.  Proposition~\ref{prop:multichain_bounded} implies that $\mathbf{\hat{0}}$ and $\mathbf{\hat{1}}$ are the least and greatest element of $\PP^{[m]}$, respectively.  

\begin{proposition}\label{prop:multichain_mobius}
	Let $\PP$ be an EL-shellable poset.  For $m>1$ we have $\mu_{\PP^{[m]}}(\mathbf{\hat{0}},\mathbf{\hat{1}})=0$.
\end{proposition}
\begin{proof}
	Theorem~\ref{thm:multichain_shellable} implies that $\PP^{[m]}$ is EL-shellable, and Theorem~\ref{thm:mobius} implies that $\mu_{\PP^{[m]}}(\mathbf{\hat{0}},\mathbf{\hat{1}})$ is essentially determined by the number of falling maximal chains (up to a sign).  Suppose that $\lambda$ is an EL-labeling of $\PP$.
	
	Let $C$ be a maximal chain in $\PP^{[m]}$.  There must be an atom $\xx$ and a coatom $\yy$ of $\PP^{[m]}$ which belong to $C$.  Corollary~\ref{cor:multichain_covers} implies that $\xx=(\hat{0},\hat{0},\ldots,\hat{0},x)$ for some atom $x$ of $\PP$, and likewise $\yy=(y,\hat{1},\hat{1},\ldots,\hat{1})$ for some coatom $y$ of $\PP$.  Since $m>1$ it follows that $\lambda^{m}(\mathbf{\hat{0}},\xx)<_{\text{lex}}\lambda(\yy,\mathbf{\hat{1}})$, which implies that $C$ cannot be falling.  
	
	We conclude that there are no falling chains in $\PP^{[m]}$, which proves the claim.
\end{proof}

\section{Applications}
	\label{sec:applications}
Let $\PP=(P,\leq)$ be a poset.  A set $I\subseteq P$ is an \defn{order ideal} of $\PP$ if for $y\in I$ and $x\leq y$ we always have $x\in I$.  Let $I(\PP)$ denote the set of order ideals of $\PP$, and let $\II(\PP)=\bigl(I(\PP),\subseteq\bigr)$.  It is quickly verified that $\II(\PP)$ is a lattice, where meet and join are given by intersection and union of sets.  

The following result together with its proof was suggested by an anonymous referee.

\begin{theorem}\label{thm:multichain_ideals}
	Let $\PP$ be a poset, and let $\CC_{m}$ be the chain with $m$ elements on the ground set $[m]$.  We have $\II(\PP\times\CC_{m})\cong\II(\PP)^{[m]}$.
\end{theorem}
\begin{proof}
	Let $\PP=(P,\leq)$.  Fix $X\in I(\PP\times\CC_{m})$, and define $X_{j}=\{i\in P\mid (i,j)\in X\}$ for $j\in[m]$.  It follows that $X_{m}\subseteq X_{m-1}\subseteq\cdots\subseteq X_{1}$, and in particular $(X_{m},X_{m-1},\ldots,X_{1})\in\II(\PP)^{[m]}$.  Let $f$ denote the map $X\mapsto (X_{m},X_{m-1},\ldots,X_{1})$.  It is straightforward to verify that $f$ is injective and order-preserving, \ie $X\subseteq Y$ implies $f(X)\subseteq f(Y)$.
	
	The inverse map $f^{-1}$ sends a multichain of order ideals $(X_{1},X_{2},\ldots,X_{m})\in\II(\PP)^{[m]}$ to
	\begin{displaymath}
		X = \bigcup_{j=1}^{m}{\bigl(X_{j}\times\{m-j+1\}\bigr)},
	\end{displaymath}
	and this map is easily seen to be order-preserving as well.
\end{proof}

In fact, Theorem~\ref{thm:multichain_ideals} tells us exactly what the poset of $m$-multichains of a distributive lattice looks like.

\begin{example}\label{ex:distributive_ideals}
	Let $\PP=(P,\leq)$ be a lattice.  An element $j\in P\setminus\{\hat{0}\}$ is \defn{join-irreducible} if whenever $j=x\vee y$, then $j\in\{x,y\}$.  Let $J(\PP)$ denote the set of join-irreducible elements of $\PP$, and let $\mathcal{J}(\PP)=\bigl(J(\PP),\leq\bigr)$.   G.~Birkhoff's representation theorem for finite distributive lattices states that $\PP$ is distributive if and only if $\PP\cong\II\bigl(\mathcal{J}(\PP)\bigr)$~\cite{birkhoff37rings}*{Theorem~5}.
	
	Theorem~\ref{thm:multichain_ideals} thus implies that if $\PP$ is a distributive lattice we have $\PP^{[m]}\cong\II\bigl(\mathcal{J}(\PP)\times\CC_{m}\bigr)$ for $m>0$.
\end{example}

The following two examples appeared as separate propositions in an earlier version of this note.

\begin{example}\label{ex:chains}
	Let $\PP=\CC_{n-1}$.  Observe that $\II(\CC_{n-1})\cong\CC_{n}$.  Theorem~\ref{thm:multichain_ideals} then implies $\CC_{n}^{[m]}\cong\II(\CC_{n-1}\times\CC_{m})$ for $n,m>0$.
\end{example}

\begin{example}\label{ex:antichains}
	Let $\PP$ be the antichain on $n$ elements, \ie the poset in which no two distinct elements are comparable.  Then $\II(\PP)=\mathcal{B}_{n}$, where $\mathcal{B}_{n}$ is the Boolean lattice of size $2^{n}$.  Moreover, $\II(\PP\times\CC_{m})\cong\CC_{m+1}^{n}$, since $\PP\times\CC_{m}$ consists of $n$ disjoint copies of $\CC_{m}$.  Any order ideal in $\PP\times\CC_{m}$ is thus composed of order ideals in each of the copies.  Theorem~\ref{thm:multichain_ideals} then implies $\mathcal{B}_{n}^{[m]}\cong\CC_{m+1}^{n}$ for $n,m>0$.
\end{example}

As pointed out by an anonymous referee, Theorem~\ref{thm:multichain_ideals} also connects the $m$-multichains of a distributive lattice to R.~Stanley's $\PP$-partitions~\cite{stanley72ordered}.  

Let $\PP=(P,\leq)$ be a finite poset with $\lvert P\rvert=p$.  A \defn{labeling} of $\PP$ is simply a bijection $\omega:P\to[p]$.  A \defn{natural labeling} is a labeling $\omega$ of $\PP$ such that $x\leq y$ implies $\omega(x)\leq\omega(y)$.  A \defn{$(\PP,\omega)$-partition} of $n$ is a map $\sigma:P\to\mathbb{N}$ such that 
\begin{enumerate}[(i)]
	\item $x\leq y$ implies $\sigma(x)\geq\sigma(y)$;
	\item $x\leq y$ and $\omega(x)>\omega(y)$ implies $\sigma(x)>\sigma(y)$;
	\item $\sum\limits_{x\in P}{\sigma(x)}=n$.
\end{enumerate}
In particular, if $\omega$ is a natural labeling, then $(\PP,\omega)$-partitions are simply order-reversing maps from $P$ to $\mathbb{N}$.  This is the case of interest here, and we will hence drop the labeling $\omega$ from the notation.  Observe further that any natural labeling produces the same set of $\PP$-partitions.

For a $\PP$-partition $\sigma$, the values $\sigma(x)$ for $x\in P$ are its \defn{parts}, and a $(\PP;m)$-partition is the a $\PP$-partition with largest part $\leq m$.  Let $A(\PP;m,n)$ be the set of all $(\PP;m)$-partitions of $n$ subject to a natural labeling of $\PP$.  Let $A(\PP;m)=\bigcup\limits_{n\geq 0}{A(\PP;m,n)}$.

Let us order the elements of $A(\PP;m)$ componentwise, \ie for $\sigma,\tau\in A(\PP;m)$ set $\sigma\preceq\tau$ if and only if $\sigma(x)\leq\tau(x)$ for all $x\in P$.  Let $\mathcal{A}(\PP;m)$ denote the resulting poset.

\begin{proposition}[\cite{stanley72ordered}*{Proposition~17.1}]\label{prop:partitions_order_ideals}
	For $\PP$ a finite poset we have $\mathcal{A}(\PP;m)\cong\II(\PP\times\CC_{m})$.
\end{proposition}

Essentially, $X\in I(\PP\times\CC_{m})$ defines $\sigma\in A(\PP;m,n)$ given by $\sigma(x)=k\geq 0$ if and only if $(x,k)\in X$ and $(x,k+1)\notin X$.  Conversely, $\sigma\in A(\PP;m,n)$ defines the order ideal $X=\bigl\{(x,k)\mid 0<k\leq\sigma(x)\bigr\}$ of $\PP\times\CC_{m}$.  

\begin{corollary}[\cite{stanley72ordered}*{Section~5}]
	For $\PP$ a finite poset we have $\mathcal{A}(\PP;m)\cong\II(\PP)^{[m]}$.
\end{corollary}

In general, a poset $\PP$ admits much less $m$-multichains than it admits $(\PP;m)$-partitions.  If we take $\PP$ to be the lattices in Figure~\ref{fig:p} and \ref{fig:q}, respectively, then Proposition~\ref{prop:partitions_order_ideals} implies that the number of $(\PP;2)$-partitions is $46$ and $33$, while the inspection of Figures~\ref{fig:p2} and \ref{fig:q2} yields that the number of $2$-multichains in these lattices is $12$ and $13$, respectively.

\section*{Acknowledgments}
I am indebted to Nathan Williams for asking whether the poset of multichains of an EL-shellable poset is again EL-shellable.  I also thank an anonymous referee for valuable comments that helped improve the exposition of this paper; in particular for suggesting Theorem~\ref{thm:multichain_ideals}.

\bibliography{../../literature}

\begin{bibdiv}
\begin{biblist}

\bib{armstrong09generalized}{article}{
      author={Armstrong, Drew},
       title={{Generalized Noncrossing Partitions and Combinatorics of Coxeter
  Groups}},
        date={2009},
     journal={Memoirs of the American Mathematical Society},
      volume={202},
}

\bib{bessis03dual}{article}{
      author={Bessis, David},
       title={{The Dual Braid Monoid}},
        date={2003},
     journal={Annales Scientifiques de l'{\'E}cole Normale Sup{\'e}rieure},
      volume={36},
       pages={647\ndash 683},
}

\bib{bessis15finite}{article}{
      author={Bessis, David},
       title={{Finite Complex Reflection Arrangements are $K(\pi,1)$}},
        date={2015},
     journal={Annals of Mathematics},
      volume={181},
       pages={809\ndash 904},
}

\bib{bessis11cyclic}{article}{
      author={Bessis, David},
      author={Reiner, Victor},
       title={{Cyclic Sieving of Noncrossing Partitions for Complex Reflection
  Groups}},
        date={2011},
     journal={Annals of Combinatorics},
      volume={15},
       pages={197\ndash 222},
}

\bib{birkhoff35on}{article}{
      author={Birkhoff, Garrett},
       title={{On the Structure of Abstract Algebras}},
        date={1935},
     journal={Proceedings of the Cambridge Philosophical Society},
      volume={31},
       pages={433\ndash 454},
}

\bib{birkhoff37rings}{article}{
      author={Birkhoff, Garrett},
       title={{Rings of Sets}},
        date={1937},
     journal={Duke Mathematical Journal},
      volume={3},
       pages={443\ndash 454},
}

\bib{bjorner80shellable}{article}{
      author={Bj{\"o}rner, Anders},
       title={{Shellable and Cohen-Macaulay Partially Ordered Sets}},
        date={1980},
     journal={Transactions of the American Mathematical Society},
      volume={260},
       pages={159\ndash 183},
}

\bib{bjorner96topological}{collection}{
      author={Bj{\"o}rner, Anders},
       title={{Topological Methods}},
      series={Handbook of Combinatorics},
   publisher={MIT Press},
     address={Cambridge, MA},
        date={1996},
      volume={2},
}

\bib{bjorner96shellable}{article}{
      author={Bj{\"o}rner, Anders},
      author={Wachs, Michelle~L.},
       title={{Shellable and Nonpure Complexes and Posets I}},
        date={1996},
     journal={Transactions of the American Mathematical Society},
      volume={348},
       pages={1299\ndash 1327},
}

\bib{bjorner97shellable}{article}{
      author={Bj{\"o}rner, Anders},
      author={Wachs, Michelle~L.},
       title={{Shellable and Nonpure Complexes and Posets II}},
        date={1997},
     journal={Transactions of the American Mathematical Society},
      volume={349},
       pages={3945\ndash 3975},
}

\bib{breckenridge93chains}{article}{
      author={Breckenridge, J.~Wylie},
      author={Varadarajan, Kalathoor},
      author={Wehrhahn, Karl~H.},
       title={{Chains, Multichains and M{\"o}bius Numbers}},
        date={1993},
     journal={Discrete Mathematics},
      volume={120},
       pages={37\ndash 50},
}

\bib{burris81course}{book}{
      author={Burris, Stanley~N.},
      author={Sankappanavar, Hanamantagouda~P.},
       title={{A Course in Universal Algebra}},
   publisher={Springer},
     address={New York},
        date={1981},
}

\bib{davey02introduction}{book}{
      author={Davey, Brian~A.},
      author={Priestley, Hilary~A.},
       title={{Introduction to Lattices and Order}},
   publisher={Cambridge University Press},
     address={Cambridge},
        date={2002},
}

\bib{gratzer11lattice}{book}{
      author={Gr{\"a}tzer, George},
       title={{Lattice Theory: Foundation}},
   publisher={Springer},
     address={Basel},
        date={2011},
}

\bib{kallipoliti15cover}{article}{
      author={Kallipoliti, Myrto},
      author={M{\"u}hle, Henri},
       title={{The $m$-Cover Posets and their Applications}},
        date={2015},
     journal={Advances in Applied Mathematics},
      volume={69},
       pages={65\ndash 108},
}

\bib{malliah86equivalence}{article}{
      author={Malliah, Chinthayamma},
      author={Bhatta, S.~Parameshwara},
       title={{Equivalence of $M$-Symmetry and Semimodularity in Lattices}},
        date={1986},
     journal={Bulletin of the London Mathematical Society},
      volume={18},
       pages={321\ndash 432},
}

\bib{rota64foundations}{article}{
      author={Rota, Gian-Carlo},
       title={{On the Foundations of Combinatorial Theory I: Theory of
  M{\"o}bius Functions}},
        date={1964},
     journal={Zeitschrift f{\"u}r Wahrscheinlichkeitstheorie und verwandte
  Gebiete},
      volume={2},
       pages={340\ndash 368},
}

\bib{stanley72ordered}{article}{
      author={Stanley, Richard~P.},
       title={{Ordered Structures and Partitions}},
        date={1972},
     journal={Memoirs of the American Mathematical Society},
      volume={119},
}

\bib{stanley74finite}{article}{
      author={Stanley, Richard~P.},
       title={{Finite Lattices and Jordan-H{\"o}lder Sets}},
        date={1974},
     journal={Algebra Universalis},
      volume={4},
       pages={361\ndash 371},
}

\bib{stanley11enumerative}{book}{
      author={Stanley, Richard~P.},
       title={{Enumerative Combinatorics, Vol. 1}},
     edition={2},
   publisher={Cambridge University Press},
     address={Cambridge},
        date={2011},
}

\bib{wachs07poset}{collection}{
      author={Wachs, Michelle~L.},
      editor={Miller, Ezra},
      editor={Reiner, Victor},
      editor={Sturmfels, Bernd},
       title={{Poset Topology: Tools and Applications}},
      series={IAS/Park City Mathematics Series},
   publisher={American Mathematical Society},
     address={Providence, RI},
        date={2007},
      volume={13},
}

\end{biblist}
\end{bibdiv}

\end{document}